\documentclass[12pt]{article}
\usepackage{amsmath,amssymb,url,color,tikz,amsthm}
\usepackage{graphicx}
\usepackage{array}
\usepackage{mathtools}
\usepackage{enumerate}
\usepackage{lscape}
\usepackage{verbatim}
\usepackage{multicol}
\usepackage{hyperref}
\newcommand{\ZZ}{\mathbb{Z}}

\newtheorem{theorem}{Theorem}[section]

\begin{document}

\begin{center}
{\bf Trimer covers in the triangular grid:} \\
{\bf twenty mostly open problems} \\
\ \\
James Propp, UMass Lowell \\
January 17, 2024
\end{center}

\noindent
\begin{abstract}
\noindent
In the past three decades, the study of rhombus tilings and domino 
tilings of various plane regions has been a thriving subfield of 
enumerative combinatorics. Physicists classify such work as the 
study of dimer covers of finite graphs. In this article we move 
beyond dimer covers to trimer covers, introducing plane regions 
called {\em benzels\/} that play a role analogous to hexagons 
for rhombus tilings and Aztec diamonds for domino tilings, inasmuch 
as one finds many (so far mostly conjectural) exact formulas 
governing the number of tilings.
\end{abstract}


\section{Introduction} \label{sec-intro}

If $G=(V,E)$ is a finite graph,
a {\bf trimer} on $G$ is a three-element subset of $V$
whose induced subgraph in $G$ is connected,
and a {\bf trimer cover} is a partition of $V$ into trimers.
Solving the trimer model for $G$
means counting the possible trimer covers.
In the physics literature one often assigns weights to the trimers,
and solving the weighted model means
finding the sum of the weights of all the trimer covers of $G$,
where the weight of an individual trimer cover
is the product of the weights of its constituent trimers;
this sum is called the {\em partition function\/}.
As a special case, one can set some of the weights equal to 1 
and the rest equal to 0; then the partition function is just
the number of trimer covers that only use ``permitted'' trimers
(trimers that have been assigned weight 1).

Physicists have studied the asymptotics of 
the weighted trimer model on the triangular lattice
and obtained formulas for the entropy in various regimes (see~\cite{VeNi}).
These results can be seen as analogous to
formulas for the entropy of dimer models on various kinds of finite graphs.
However, as far as I am aware there have not been
(up to now) any exact enumerative results
for trimer models on finite subgraphs of the triangular lattice
in the style of the well-known enumerations of 
rhombus tilings of hexagons and domino tilings of the Aztec diamond~\cite{Pro1}.

\begin{figure}[h]
\begin{center}
\includegraphics[width=3.0in]{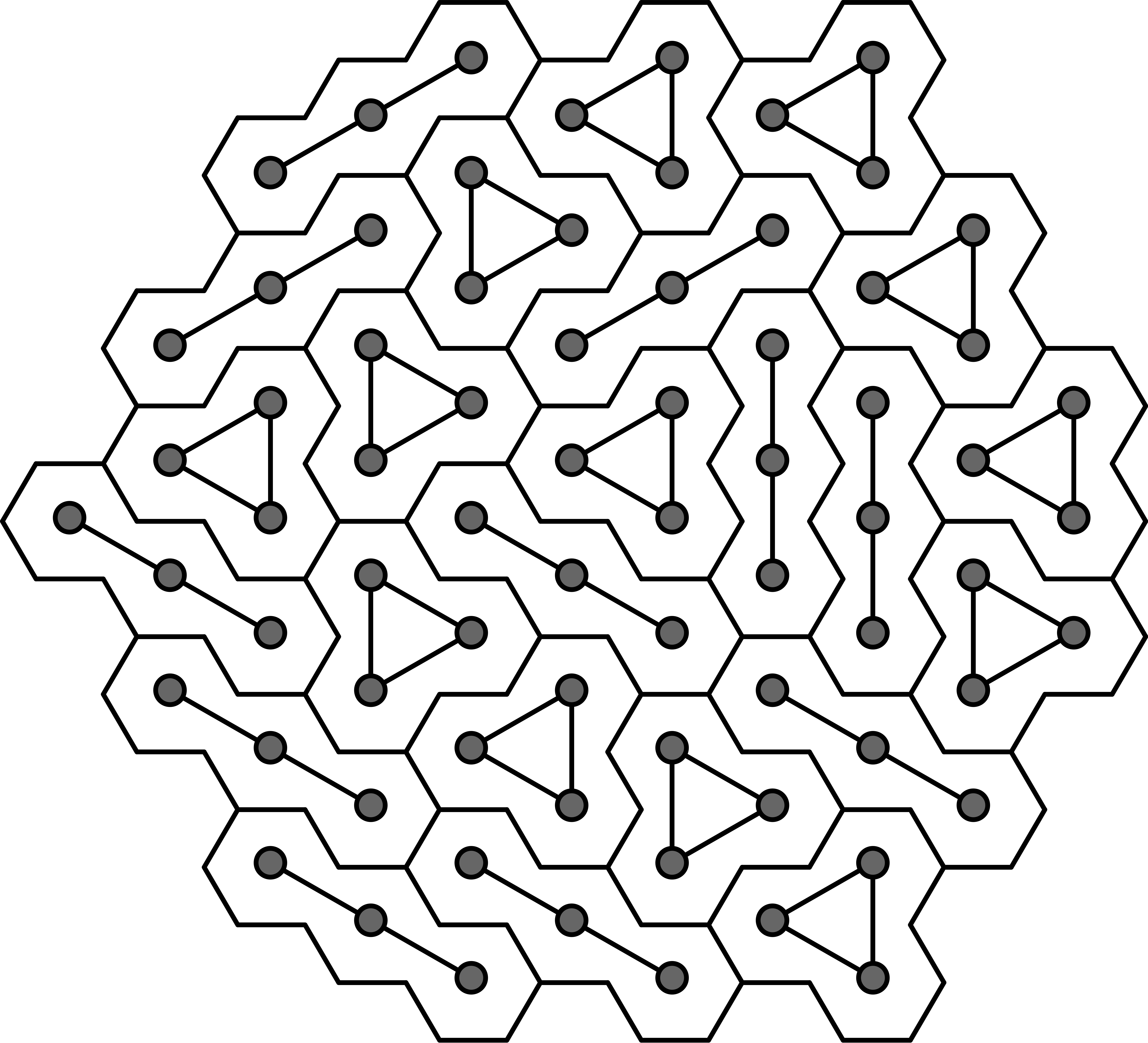}
\end{center}
\caption{The (9,9)-benzel tiled by stones and bones
and the associated trimer cover of the (9,9)-benzel graph.}
\label{fig:sample-tiling}
\end{figure}

In this article I introduce 
finite subgraphs of the triangular lattice 
that should interest enumerative combinatorialists
inasmuch as the number of trimer covers
appears to be given by exact formulas in many cases.
See Figure~\ref{fig:sample-tiling} 
which shows both a trimer cover of a finite subgraph
of the triangular lattice
and the associated tiling of a dual 2-complex
whose faces correspond to vertices of the triangular lattice.
In tiling language, the tiled region in the Figure
is a {\em polyhex} and the tiles are {\em trihexes}.
I have dubbed the tiled region a {\em benzel}, 
and the allowed tiles {\em stones} and {\em bones}.
I will switch between the trimer picture
and the tiling picture throughout the article.
The graph admits an essentially unique coloring
with the colors red, blue, and green
in which no two adjacent vertices have been assigned the same color;
under such a coloring, each stone-trimer and bone-trimer
contains exactly one vertex of each color.

This article is a companion to the talk I gave at
the Open Problems in Algebraic Combinatorics conference 
on May 18, 2022; the slides and video 
are available through the conference website
\vspace{-10pt}
\begin{center}
\href{http://www.samuelfhopkins.com/OPAC/opac.html}{{\tt http://www.samuelfhopkins.com/OPAC/opac.html}}
\end{center}
\vspace{-10pt}
and updates on the problems are available at
\vspace{-10pt}
\begin{center}
\href{http://jamespropp.org/benzels.html}{{\tt http://jamespropp.org/benzels.html}}.
\end{center}
\vspace{-10pt}

\begin{figure}[h]
\begin{center}
\includegraphics[width=5.6in]{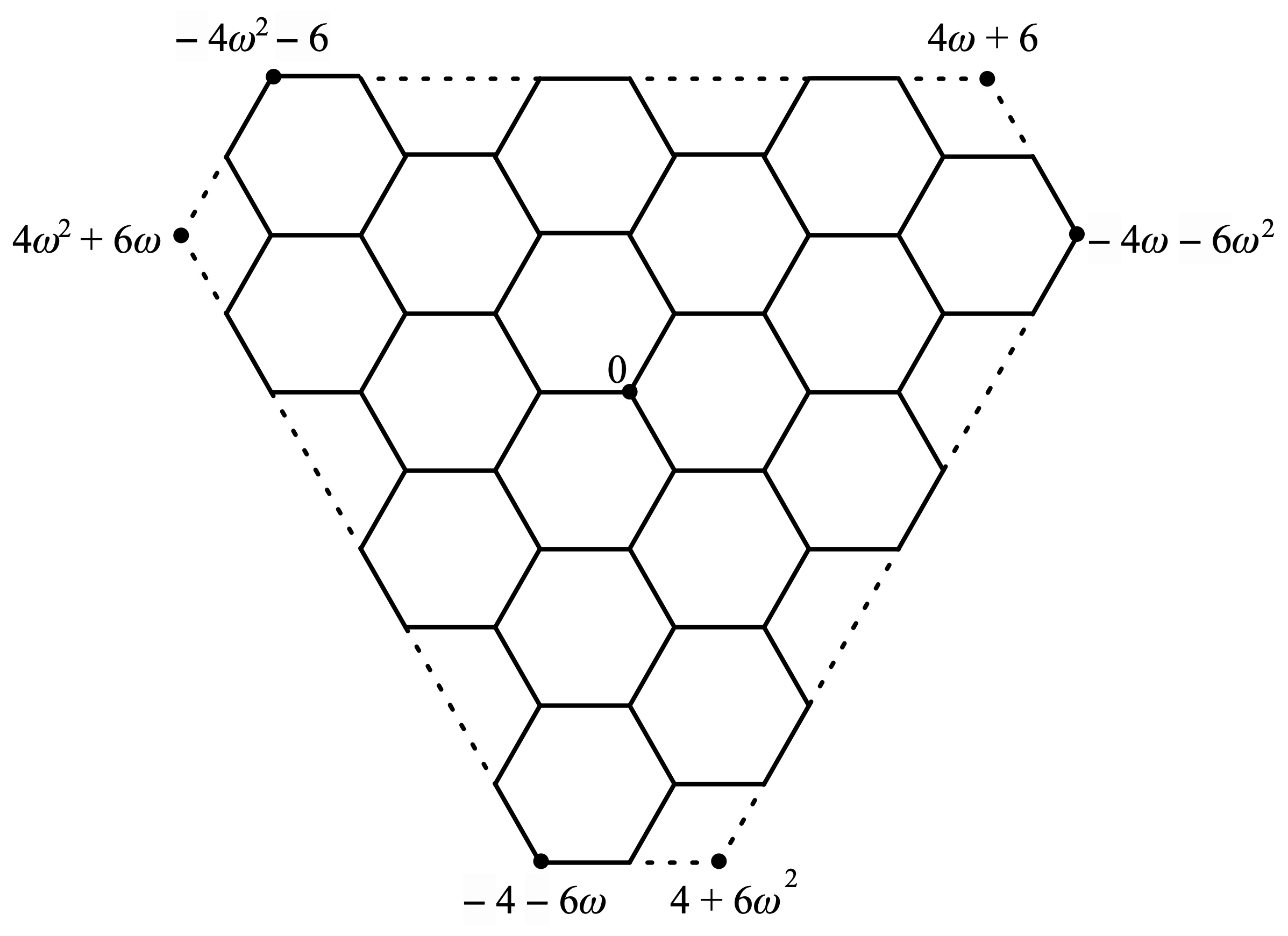}
\end{center}
\caption{The (4,6)-benzel and its enclosing hexagon.}
\label{fig:one-benzel}
\end{figure}

Consider the complex plane tiled by unit hexagonal cells
centered at 1, $\omega$, and $\omega^2$ (here and hereafter 
$\omega$ denotes $e^{2 \pi i / 3}$);
the cell centered at $\alpha$ has corners at 
$\alpha \pm 1$, $\alpha \pm \omega$, and $\alpha \pm \omega^2$.
Given positive integers $a,b$ satisfying
$2 \leq a \leq 2b$ and $2 \leq b \leq 2a$, we define 
the $\boldsymbol{(a,b)}$-{\bf benzel} as the union of the cells
that lie fully inside the hexagon with vertices
$a\omega+b$, $-a\omega^2-b$, $a\omega^2+b\omega$, $-a-b\omega$,
$a+b\omega^2$, and $-a\omega-b\omega^2$
(a hexagon centered at 0 with threefold rotational symmetry
whose side-lengths alternate between $2a-b$ and $2b-a$,
degenerating to a triangle when $a=2b$ or $b=2a$),
as shown in Figure~\ref{fig:one-benzel} for $a=4$, $b=6$.
(Mnemonically, $a$ is the height of the top edge of the hexagon 
{\em above} 0, and $b$ is the depth of the bottom edge of the hexagon 
{\em below} 0, using appropriately scaled units.) We lose no generality 
in assuming $2 \leq a \leq 2b$ and $2 \leq b \leq 2a$,
since the $(a,b)$-benzel as defined above
coincides with the $(a,a-b)$-benzel when $a > 2b$
and with the $(b-a,b)$-benzel when $b > 2a$,
the former satisfying $a \leq 2(a-b)$
and the latter satisfying $b \leq 2(b-a)$.

\begin{figure}[h]
\begin{center}
\includegraphics[width=4.0in]{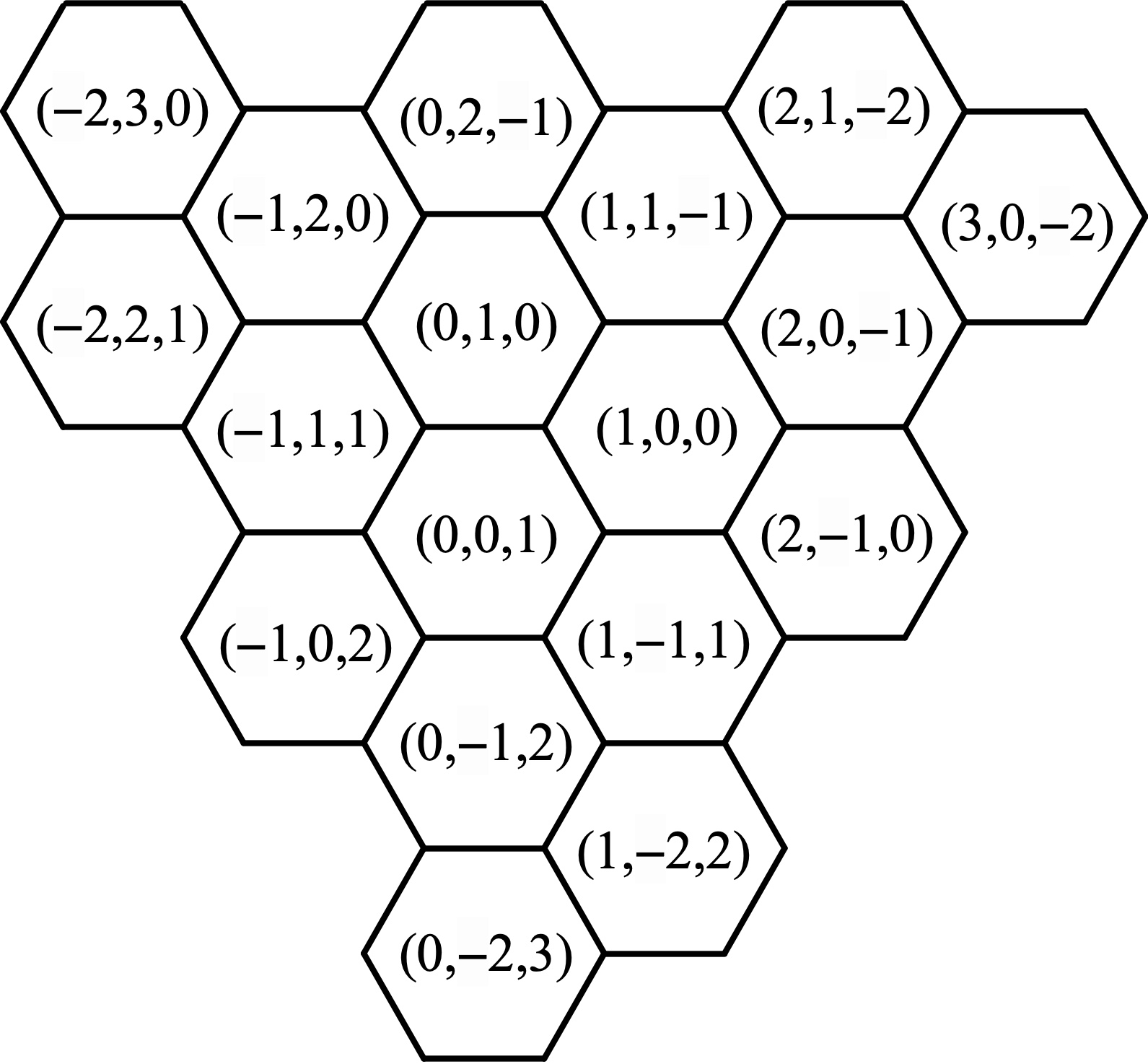}
\end{center}
\caption{Barycentric coordinates for the cells of the (4,6)-benzel.}
\label{fig:barycentric}
\end{figure}

Here is an alternative description of benzels
in terms of the centers of the cells using barycentric coordinates
relative to the triangle with vertices 1, $\omega$, and $\omega^2$.
Each cell's center point $\alpha$ belongs to $\ZZ[\omega]$
and can be represented by the unique $(i,j,k) \in \ZZ^3$ 
satisfying $i+j\omega+k\omega^2=\alpha$ and $i+j+k=1$
(here and hereafter $i$ denotes a nonnegative integer
as opposed to the square root of $-1$).
Then the $(a,b)$-benzel consists of those cells
whose centers $(i,j,k)$ satisfy
$-(a\!-\!1) \leq j\!-\!i,k\!-\!j,i\!-\!k \leq b\!-\!1$.
Figure~\ref{fig:barycentric} shows the (4,6)-benzel with its cells marked 
with the barycentric coordinates of their respective center points.

Let $V = \{(i,j,k) \in \ZZ^3 \: : \: i+j+k=1, 
\ -(a\!-\!1) \leq j\!-\!i,\,k\!-\!j,\,i\!-\!k \leq b\!-\!1\}$.
Given $(i_1,j_1,k_1)$ and $(i_2,j_2,k_2)$ in $V$,
join $(i_1,j_1,k_1)$ and $(i_2,j_2,k_2)$ by an edge
when $|i_1\!-\!i_2|+|j_1\!-\!j_2|+|k_1\!-\!k_2|=2$
(that is, when the unit hexagons centered on those two vertices share an edge).
This is the $\boldsymbol{(a,b)}$-{\bf benzel graph}.

The $(a,b)$-benzel has threefold rotational symmetry
but for most $a,b$ it does not have bilateral symmetry.
Exchanging $a$ and $b$ corresponds to reflecting the benzel
across a horizontal axis (or if one prefers across an axis
making a 60 degree angle with the horizontal axis).

\begin{figure}[h]
\begin{center}
\includegraphics[width=5.0in]{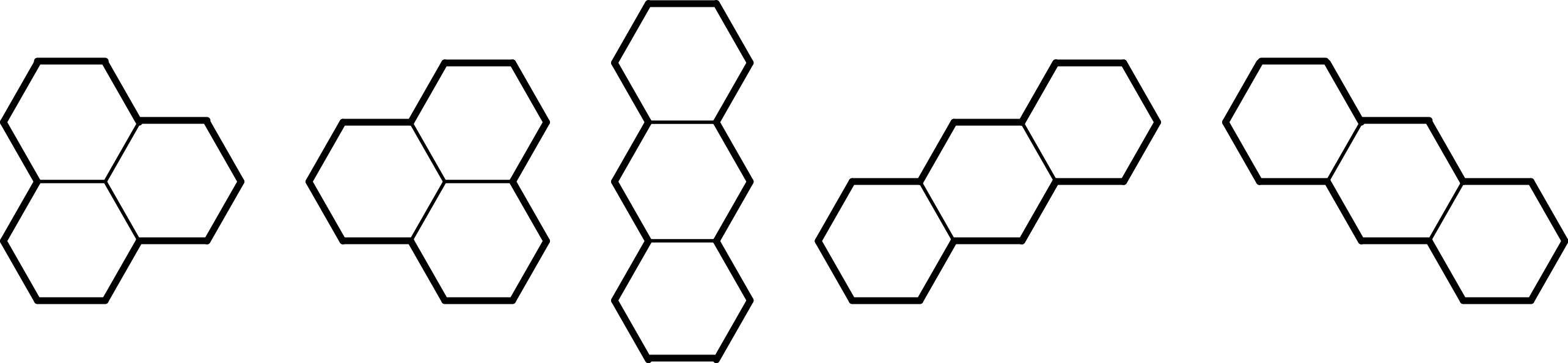}
\end{center}
\caption{The five prototiles: the right stone, the left stone,
the vertical bone, the rising bone, and the falling bone.}
\label{fig:prototiles}
\end{figure}

We consider tilings of the $(a,b)$-benzel
by way of five sorts of prototiles.
These prototiles (shown in Figure~\ref{fig:prototiles}) are 
the {\bf right(-pointing) stone}, the {\bf left(-pointing) stone},
the {\bf vertical bone}, the {\bf rising bone}, and the {\bf falling bone}.
The right stone is a benzel
(specifically, the (2,2)-benzel) but the left stone is not.
As usual, the allowed tiles are the translates of the prototiles. 
Dually we form spanning subgraphs of the $(a,b)$-benzel graph
whose connected components all consist of three vertices.
Figure~\ref{fig:sample-tiling} shows a tiling of the (9,9)-benzel 
and the associated trimer cover of the (9,9)-benzel graph.

There is a third kind of trihex
that consists of three hexagons whose centers form a 120-degree angle. 
Timothy Chow has dubbed this kind of tile a {\em phone}.
Since the associated trimers do not respect 
the tripartite nature of the triangular lattice,
I have chosen to exclude those trimers from consideration here.
This is not to suggest that such trimers are unworthy of study.
Just as enumeration of dimer covers 
should not be limited to the bipartite case,
enumeration of trimer covers 
should not be limited to the tripartite case.
However, it seems best to start with the problems
that are likely to be easier.
In particular, the Conway-Lagarias invariant discussed below
only applies if all our trihexes are stones and bones.

Conway and Lagarias~\cite{CoLa} studied tilings like the ones considered here
(calling stones and bones $T_2$ and $L_3$ tiles respectively).
They showed that for any simply-connected region
in the hexagonal grid that can be tiled by stones and bones,
the total area of the right stones minus the total area of the left stones 
does not depend on the specific tiling
but only depends on the region being tiled.
This is the {\em Conway-Lagarias invariant} of the region.
It can be shown (see the companion article~\cite{KiPr})
that the area of the $(a,b)$-benzel (as measured by the number of tiles) 
is given by
$$
\begin{array}{ll}
(- a^2 + 4 a b - b^2 - a - b)/2 & \mbox{if $a+b \equiv 0$ or 2 (mod 3)}, \\
(- a^2 + 4 a b - b^2 - a  - b + 2)/2 & \mbox{if $a+b \equiv 1$ (mod 3)} \\
\end{array}
$$
while the value of the Conway-Lagarias invariant is
$$
\begin{array}{ll}
(3a^2 - 6ab + 3b^2 - a - b)/2
& \mbox{if $a+b \equiv 0$ (mod 3)}, \\
(-a^2 + 4ab - b^2 - a - b + 2)/2
& \mbox{if $a+b \equiv 1$ (mod 3)}, \\
(3a^2 - 6ab + 3b^2 + a + b - 2)/2
& \mbox{if $a+b \equiv 2$ (mod 3).} \\
\end{array} 
$$
(The expressions in Theorem 2 of~\cite{KiPr} have opposite sign
because the benzels in that article are mirror images
of the benzels considered here.)
Note that when $a+b \equiv 1$ (mod 3),
the Conway-Lagarias invariant is equal to the area of the benzel
so that the tiling must consist entirely of right-pointing stones;
for instance, this is the case with the (4,6)-benzel shown earlier.

The set of 5 prototiles has $2^5-1 = 31$ nonempty subsets,
and for each, we can ask 
in how many ways it is possible to tile the $(a,b)$-benzel,
that is, in how many ways it is possible 
find translates of the prototiles
whose interiors are disjoint and whose union is the benzel.
There is some redundancy here.
Because the benzel has threefold rotational symmetry,
and because 120 degree rotations preserve the two stones' orientation
(right versus left),
the number of tilings depends only on
(a) whether right stones are allowed,
(b) whether left stones are allowed,
and (c) how many of the three kinds of bones are allowed (0, 1, 2, or 3).
Thus there are really only $(2)(2)(4)-1 = 15$ 
classes of tiling problems to consider.
For $0 \leq i,j \leq 1$ and $0 \leq k \leq 3$
we define $T_{ijk}(a,b)$ as the number of ways to tile the $(a,b)$-benzel
if the set of allowed prototiles contains the right stone iff $i=1$,
contains the left stone iff $j=1$, and contains $k$ of the bones.
We say that such a tiling is an $i,j,k$ tiling.

It is not hard to show 
that for each of the 15 cases, $T_{ijk}(a,b) = T_{ijk}(b,a)$.
It is also not hard to show that the $(n,2n)$-benzel is the same 
as the $(n,2n-1)$-benzel and the $(n,2n-2)$-benzel.
Consequently, in the tables that follow 
I provide values for $T_{ijk}(a,b)$
only for $2 \leq a \leq 2b-2$ and $2 \leq b \leq 2a-2$. 
Lastly, it is not hard to show that 
when $i=0$ and $a+b \equiv 2$ (mod 3), $T_{ijk}(a,b) = 0$.
That is because in this case the Conway-Lagarias invariant
is equal to the positive number $(3(a-b)^2 + (a+b-2))/2$,
implying that every tiling of the $(a,b)$-benzel
must have at least one right-pointing stone.

David desJardins wrote a general purpose program {\tt TilingCount}
that I used to enumerate tilings of regions
with various sets of allowed prototiles.
This led to the questions and conjectures that appear below.
I am happy to share the code and the data on which my conjectures are based
(much of which can be found in the Online Encyclopedia of Integer Sequences).
Table~\ref{tab:map} is a map of 
the first eighteen problems presented in this article
as they relate to those fifteen cases.
Rows describe which stones are allowed;
columns describe how many bones are allowed.
(This table omits cases where 
the number of allowed bone prototiles is zero or one;
in such situations at most one tiling exists, 
even when both of the stone prototiles are allowed.)

\begin{table} 
\begin{center}
\caption{A map of the first eighteen problems.}
\label{tab:map}
\vspace{0.2in}
\begin{tabular}{c|c|c}
\                              & \rm Two \ kinds \ of \ bones    & \rm Three \ kinds \ of \ bones \\
\hline
\rm No \ stones                & \mbox{(no tilings exist)}       & \mbox{type 003: prob.\ 1}     \\
\hline
\rm Left \ stones          & \mbox{type 012: probs.\ 2--3}  & \mbox{type 013: prob.\ 4}     \\
\hline
\rm Right \ stones          & \mbox{type 102: prob.\ 5}      & \mbox{type 103: probs.\ 6--7}  \\
\hline
\rm Both \ kinds \ of \ stones & \mbox{type 112: probs.\ 8--13} & \mbox{type 113: probs.\ 14--18}
\end{tabular}
\end{center}
\end{table}

Benzels behave differently according to whether $a+b$ is 0, 1, or 2 (mod 3),
so in what follows we will sometimes subdivide analyses
according to the congruence class of $a+b$.

In the cases where only two kinds of bone tiles are permitted,
the allowed tilings can be viewed as ribbon tilings, as in~\cite{Pak};
indeed, this was the mode of presentation employed in~\cite{CoLa},
with four of the five prototiles being depicted as ribbons.
Switching over to the ribbon tilings presentation
has already yielded solutions to problems 2 and 3 in~\cite{DLPY},
and is likely to provide leverage 
on other problems in the first column of the table.

\section{No stones, three kinds of bones} \label{sec-no-three}

Prior to the conference, I was able to show that if 
an $(a,b)$-benzel can be tiled by bones alone, then we must have 
$a = k(3k-1)/2$ and $b = k(3k+1)/2$ (or vice versa) for some $k \geq 2$.
Several attending students found a proof
that this necessary condition is also sufficient.
(Note that such benzels belongs to the case $a+b \equiv 0$ (mod 3).)
Jesse Kim found the most complete solution, providing an explicit proof
that the tiling he described works for all $k$.
It appears (see \href{https://oeis.org/A364134}{OEIS entry A364134})
that the number of such tilings grows exponentially in $k^4$.


{\bf Problem 1:} Find an exact formula for $T_{003}(k(3k-1)/2,k(3k+1)/2)$.


Of course, even short of an exact formula,
any method of determining the number of tilings
that is more efficient than brute-force enumeration
(e.g., a recurrence relation) would be of interest.

See~\cite{KiPr} for more discussion of no-stones tilings of benzels.

\section{Left stones, two kinds of bones} \label{sec-neg-two}

Table~\ref{tab:012} shows the values of $T_{012}(a,b)$ for $a,b \leq 10$.

\begin{table} 
\begin{center}
\caption{Values of $T_{012}(a,b)$.}
\label{tab:012}
\vspace{0.2in}
\begin{tabular}{c|ccccccccc}
$a \backslash b$  & 2 & 3 & 4 & 5 & 6 & 7 & 8 & 9 & 10 \\
\hline
 2       & 0 &   &   &   &   &   &   &   &    \\
 3       &   & {\bf 2} & 0 &   &   &   &   &   &    \\
 4       &   & 0 & 0 & {\bf 2} & 0 &   &   &   &    \\
 5       &   &   & {\bf 2} & 0 & 0 & 0 & 0 &   &    \\
 6       &   &   & 0 & 0 & {\bf 8} & 0 & 0 & 0 & 0  \\
 7       &   &   &   & 0 & 0 & 0 & {\bf 8} & 0 & 0  \\
 8       &   &   &   & 0 & 0 & {\bf 8} & 0 & 0 & 0  \\
 9       &   &   &   &   & 0 & 0 & 0 & {\bf 48} & 0  \\
10       &   &   &   &   & 0 & 0 & 0 & 0 & 0  
\end{tabular}
\end{center}
\end{table}


{\bf Problem 2:} Is it true that $T_{012}(3n,3n) = 2^n n!$ for $n \geq 1$?


{\bf Problem 3:} Is it true that $T_{012}(3n+1,3n+2) = 2^n n!$ for $n \geq 1$?


Comment: Colin Defant, Rupert Li, Benjamin Young and I 
worked on problem 2 during the conference 
and later succeeded in solving Problems 2 and 3; see~\cite{DLPY}.
We were also able to verify that the pattern of entries equal to 0 and 1
in Tables~\ref{tab:012} through~\ref{tab:102}
persists for all larger values of $a$ and $b$.

\section{Left stones, three kinds of bones} \label{sec-neg-three}

Table~\ref{tab:013} shows the values of $T_{013}(a,b)$ for $a,b \leq 10$.

\begin{table} 
\begin{center}
\caption{Values of $T_{013}(a,b)$.}
\label{tab:013}
\vspace{0.2in}
\begin{tabular}{c|ccccccccc}
$a \backslash b$ & 2 & 3 & 4 & 5 &   6 &    7 &    8 &     9 & 10  \\
\hline
 2             & 0 &   &   &   &     &      &      &       &     \\
 3             &   & {\bf 3} & 0 &   &     &      &      &       &     \\
 4             &   & 0 & 0 & {\bf 9} &   0 &      &      &       &     \\
 5             &   &   & {\bf 9} & 0 &   0 &    {\bf 2} &    0 &       &     \\
 6             &   &   & 0 & 0 & {\bf 144} &    0 &    0 &     0 &   0 \\
 7             &   &   &   & {\bf 2} &   0 &    0 & {\bf 1143} &     0 &   0 \\
 8             &   &   &   & 0 &   0 & {\bf 1143} &    0 &     0 & {\bf 825} \\
 9             &   &   &   &   &   0 &    0 &    0 & {\bf 73454} &   0 \\
10             &   &   &   &   &   0 &    0 &  {\bf 825} &     0 &   0 
\end{tabular}
\end{center}
\end{table}

When $a$ and $b$ are such that 
the Conway-Lagarias invariant is strictly positive,
the $(a,b)$-benzel cannot be tiled by bones and left stones;
the corresponding entries in the table must be zero.
On the other hand, when $a$ and $b$ are such that 
the Conway-Lagarias invariant is negative or zero,
the entries in the table are observed to be positive,
though I see no reason for concluding that they are.

{\bf Problem 4:} Is it true that when the Conway-Lagarias invariant
associated with the $(a,b)$-benzel is negative or zero,
tilings of type 013 exist?

\section{Right stones, two kinds of bones} \label{sec-pos-two}

Table~\ref{tab:102} shows the values of $T_{102}(a,b)$ for $a,b \leq 10$.

\begin{table} 
\begin{center}
\caption{Values of $T_{102}(a,b)$.}
\label{tab:102}
\vspace{0.2in}
\begin{tabular}{c|ccccccccc}
$a \backslash b$ & 2 & 3 & 4 & 5 & 6 & 7 & 8 & 9 & 10 \\
\hline
 2             & 1 & &   &   &   &   &   &   &    \\
 3             &   & 0 & 1 & &   &   &   &   &    \\
 4             &   & 1 &{\bf 2}& 0 & 1 & &   &   &    \\
 5             &   &   & 0 & 1 &{\bf 4} & 0 & 1 & &   \\
 6             &   &   & 1 & 4 & 0 & 1 &{\bf 10}& 0 & 1  \\
 7             &   &   &   & 0 & 1 &{\bf 8}& 0 & 1 &\bf{28} \\
 8             &   &   &   & 1 &10 & 0 & 1 &{\bf 24}& 0  \\
 9             &   &   &   &   & 0 & 1 &24 & 0 & 1  \\
10             &   &   &   &   & 1 &28 & 0 & 1 &{\bf 48}
\end{tabular}
\end{center}
\end{table}


{\bf Problem 5:} Is it true that 
$$T_{102}(n+3k,2n+3k-1) = 
\prod_{i=1}^k \frac{(2i)! (2i+2n-2))!}{(i+n-1)! (i+n+k-1)!}$$
for $k \geq 0$ and $n \geq 1$ (except $(k,n) = (0,1))$? 
Equality has been verified for $0 \leq k \leq 5$, $1 \leq n \leq 5$.


Comment: This formula and the $a,b$ symmetry relation
together provide a conjectural enumeration 
of tilings of type 102 of the $(a,b)$-benzel
for all $a,b$ satisfying $a+b \equiv 2$ (mod 3).
(When $a+b \equiv 1$, the number of tilings is 0;
when $a+b \equiv 2$, the number of tilings is 1.)

Comment: Three special cases merit special attention.
When $k=1$, the right-hand side of the equation
is twice the $n$th Catalan number;
when $n=1$, the right-hand side of the equation is $2^k k!$;
and when $n=2$, the right-hand side of the equation is $2^k (k+1)!$.

\section{Right stones, three kinds of bones} \label{sec-pos-three}

Table~\ref{tab:103} shows the values of $T_{103}(a,b)$ for $a,b \leq 10$.

\begin{table} 
\begin{center}
\caption{Values of $T_{103}(a,b)$.}
\label{tab:103}
\vspace{0.2in}
\begin{tabular}{c|ccccccccc}
$a \backslash b$ & 2 & 3 & 4 &  5 &   6 &   7 &     8 &     9 &     10 \\
\hline
 2             & 1 &   &   &    &     &     &       &       &        \\
 3             &   &{\bf 0}& 1 &    &     &     &       &       &        \\
 4             &   & 1 & 7 &{\bf 0}&  1  &     &       &       &        \\
 5             &   &   & 0 &  1 & 33  &{\bf 2} &     1 &       &        \\
 6             &   &   & 1 & 33 &  0  &  1  &   164 &   {\bf 21}&      1 \\
 7             &   &   &   &  2 &  1  & 666 &     0 &     1 &    864 \\
 8             &   &   &   &  1 & 164 &  0  &     1 & 12430 &      0 \\
 9             &   &   &   &    & 21  &  1  & 12430 &     0 &      1 \\
10             &   &   &   &    &  1  & 864 &     0 &     1 & 655721
\end{tabular}
\end{center}
\end{table}

{\bf Problem 6:} Is is true that 
$T_{103}(n,2n-3) = (3n+3)(3n-7)!/(n-5)!(2n-1)!$ for $n \geq 5$?
(The formula also works for $n=3$ and $n=4$
if one treats $1/(-1)!$ and $1/(-2)!$ as 0.)
Equality has been verified for $5 \leq n \leq 16$.


Aside from the fact that Problems 2 and 6 involve different prototile sets,
the two problems differ in another important way:
in problem 2 we have $b/a \rightarrow 1$ as $n \rightarrow \infty$
while in problem 6 we have $b/a \rightarrow 2$ as $n \rightarrow \infty$.
In the former case we say that the sequence is associated with 
a {\bf central diagonal} of the table of values 
while in the latter case we say that the sequence 
is associated with a {\bf peripheral diagonal}
(recall that $a/b$ and $b/a$ cannot exceed 2).


In parallel with the observations that preceded Problem 4,
note that when the Conway-Lagarias invariant of a benzel is strictly negative,
the benzel cannot be tiled by bones and right-pointing stones;
the corresponding entries in the table must be zero.
On the other hand, when $a$ and $b$ are such that 
the Conway-Lagarias invariant is positive or zero,
the entries in the table are observed to be positive,
though I see no reason for concluding that they are.

{\bf Problem 7:} Is it true that when the Conway-Lagarias invariant
associated with the $(a,b)$-benzel is positive or zero,
tilings of type 103 exist?


\section{Both stones, two kinds of bones} \label{sec-both-two}

Table~\ref{tab:112} shows the values of $T_{112}(a,b)$ for $a,b \leq 10$.
Table~\ref{tab:3n3n} gives the first few values 
of $T_{112}(3n,3n)$ in factored form.

\begin{table} 
\begin{center}
\caption{Values of $T_{112}(a,b)$.}
\label{tab:112}
\vspace{0.2in}
\begin{tabular}{c|ccccccccc}
$a \backslash b$ & 2 & 3 & 4 &  5 &  6 &   7 &    8 &     9 &     10  \\
\hline
 2             & 1 & &   &    &    &     &      &       &         \\
 3             & &{\bf 2}& 1 &    &    &     &      &       &         \\
 4             &   & 1 & {\bf 4} & {\bf 6} &  1 &     &      &       &         \\
 5             &   &   & 6 &  1 & {\bf 16} & {\bf 22} &    1 &       &         \\
 6             &   &   & 1 & 16 & {\bf 48} &   1 & {\bf 68} & {\bf 90}&      1  \\
 7             &   &   &   & 22 &  1 & {\bf 224} & {\bf 512} &     1 & {\bf 304} \\
 8             &   &   &   &  1 & 68 & 512 &    1 & {\bf 3360} &   6736  \\
 9             &   &   &   &    & 90 &   1 & 3360 & {\bf 15360} &      1  \\
10             &   &   &   &    &  1 & 304 & 6736 &     1 & {\bf 168960} 
\end{tabular}
\end{center}
\end{table}

\begin{table} 
\begin{center}
\caption{Values of $T_{112}(3n,3n)$.}
\label{tab:3n3n}
\vspace{0.2in}
\begin{tabular}{c}
$2^{1},$ \\
$2^{4} \: 3^{1},$ \\ 
$2^{10} \: 3^{1} \: 5^{1},$ \\
$2^{16} \: 7^{1} \: 11^{1} \: 13^{1},$ \\ 
$2^{28} \: 3^{2} \: 7^{1} \: 13^{1} \: 17^{1},$ \\
$2^{38} \: 3^{2} \: 11^{1} \: 17^{2} \: 19^{2},$ \\
$2^{50} \: 5^{1} \: 11^{2} \: 13^{1} \: 17^{1} \: 19^{2} \: 23^{2},$ \\ 
$2^{64} \: 3^{3} \: 5^{4} \: 11^{1} \: 13^{2} \: 19^{1} \: 23^{3} \: 29^{1},$ \\
$2^{84} \: 3^{4} \: 5^{3} \: 13^{2} \: 17^{1} \: 23^{2} \: 29^{3} \: 31^{2},$ \dots
\end{tabular}
\end{center}
\end{table}


{\bf Problem 8:} 
With $T(n)$ denoting $T_{112}(3n,3n)$,
is it true that the second quotient $T(n)T(n+2)/T(n+1)^2$ is equal to
$$\frac{256 (2n+3)^2 (4n+1) (4n+3)^2 (4n+5)}{27 (3n+1) (3n+2)^2 (3n+4)^2 (3n+5)}$$
for all $n \geq 1$? 
Equality has been verified for $1 \leq n \leq 7$. See 
\href{https://oeis.org/A352207}{OEIS entry A352207}.

Comment: David desJardins found the pattern governing the numbers $T_{112}(n)$,
with assistance from Christian Krattenthaler, Greg Kuperberg
and other members of the {\tt domino} listserv.
The same is true for Problem 10.


Table~\ref{tab:3np13np1} gives the first few values 
of $T_{112}(3n+1,3n+1)$ in factored form.

\begin{table} 
\begin{center}
\caption{Values of $T_{112}(3n+1,3n+1)$.}
\label{tab:3np13np1}
\vspace{0.2in}
\begin{tabular}{c}
$2^{2},$ \\
$2^{5} \: 7^{1},$ \\
$2^{10} \: 3^{1} \: 5^{1} \: 11^{1},$ \\
$2^{17} \: 7^{1} \: 11^{1} \: 13^{2},$ \\
$2^{30} \: 3^{1} \: 13^{1} \: 17^{2} \: 19^{1},$ \\
$2^{38} \: 3^{1} \: 11^{1} \: 17^{2} \: 19^{3} \: 23^{1},$ \dots
\end{tabular}
\end{center}
\end{table}

With $T(n) = T_{112}(3n+1,3n+1)$ with $n \geq 1$,
it appears that $T(n)$ has no prime factor greater than or equal to $4n$.

{\bf Problem 9:} 
Find a formula governing $T_{112}(3n+1,3n+1)$. See 
\href{https://oeis.org/A364481}{OEIS entry A364481}.


Table~\ref{tab:3np13np2} gives the first few values 
of $T_{112}(3n+1,3n+2)$ in factored form.

\begin{table} 
\begin{center}
\caption{Values of $T_{112}(3n+1,3n+2)$.}
\label{tab:3np13np2}
\vspace{0.2in}
\begin{tabular}{c}
$2^{1} \: 3^{1},$ \\
$2^{9},$ \: \\
$2^{9} \: 3^{1} \: 5^{1}\: 7^{1}\: 11^{1},$ \\
$2^{25} \: 3^{1} \: 7^{1} \: 13^{1},$ \\
$2^{28} \: 3^{2} \: 11^{1} \: 13^{1} \: 17^{2} \: 19^{1},$ \\
$2^{50} \: 3^{1} \: 11^{1} \: 17^{1} \: 19^{2} \: 23^{1},$ \\
$2^{49} \: 3^{2} \: 5^{2} \: 11^{2} \: 13^{2} \: 17^{1} \: 19^{2} \: 23^{3},$ \\
$2^{81} \: 3^{3} \: 5^{4} \: 13^{2} \: 23^{2} \: 29^{2} \: 31^{1},$ \dots
\end{tabular}
\end{center}
\end{table}

{\bf Problem 10:} 
With $T(n)$ denoting $T_{112}(3n+1,3n+2)$,
is it true that $T(n)T(n+3)/T(n+1)T(n+2)$ is equal to
$$\frac{65536 (2n+3)(2n+5)^2 (2n+7) (4n+3) (4n+5)^2 (4n+7)^2 (4n+9)^2 (4n+11)}{729 
(3n+2)(3n+4)^3 (3n+5)^2 (3n+7)^2 (3n+8)^3 (3n+10)}$$
for all $n \geq 1$?
Equality has been verified for $1 \leq n \leq 8$.  See 
\href{https://oeis.org/A364482}{OEIS entry A364482}.


Table~\ref{tab:3nm13n} gives the first few values 
of $T_{112}(3n-1,3n)$ in factored form.

\begin{table} 
\begin{center}
\caption{Values of $T_{112}(3n-1,3n)$.}
\label{tab:3nm13n}
\vspace{0.2in}
\begin{tabular}{c}
$1^{1},$ \\
$2^{4},$ \\
$2^{5} \: 3^{1} \: 5^{1} \: 7^{1},$ \\
$2^{16} \: 11^{1} \: 13^{1},$ \\
$2^{19} \: 3^{1} \: 7^{1} \: 11^{1} \: 13^{2} \: 17^{1},$ \\
$2^{39} \: 3^{1} \: 17^{2} \: 19^{2},$ \\
$2^{37} \: 5^{1} \: 11^{2} \: 13^{1} \: 17^{2} \: 19^{3} \: 23^{2},$ \dots
\end{tabular}
\end{center}
\end{table}

With $T(n) = T_{112}(3n-1,3n)$ with $n \geq 1$,
it appears that $T(n)$ has no prime factor greater than or equal to $4n$.

{\bf Problem 11:} 
Find a formula governing $T_{112}(3n-1,3n)$. See 
\href{https://oeis.org/A364483}{OEIS entry A364483}.


We now switch from central diagonals to peripheral diagonals.


{\bf Problem 12:} 
Is it true that $T_{112}(n+2,2n+1)$ 
is the ``$n$th large Schr\"oder number'' 
(see \href{https://oeis.org/A006318}{OEIS entry A006318})
for all $n \geq 1$?
Equality has been verified for $1 \leq n \leq 15$.


{\bf Problem 13:} 
Is it true that $T_{112}(n+2,2n)$ 
is the number of ``royal paths in a lattice of order $n$''
(see \href{https://oeis.org/A006319}{OEIS entry A006319}) 
for all $n \geq 1$?
Equality has been verified for $1 \leq n \leq 15$.

\section{Both stones, three kinds of bones} \label{sec-both-three}

Finally we come to the most permissive situation: all prototiles are allowed.
Table~\ref{tab:113} shows the values of $T_{113}(a,b)$ for $a,b \leq 10$.

\begin{table} 
\begin{center}
\caption{Values of $T_{113}(a,b)$.}
\label{tab:113}
\vspace{0.2in}
\begin{tabular}{c|ccccccccc}
$a \backslash b$ & 2 & 3 &  4 &   5 &    6 &     7 &       8 & 9       & 10      \\
\hline
 2             & 1 &   &    &     &      &       &         &         &            \\
 3             &   & 3 &  1 &     &      &       &         &         &            \\
 4             &   & 1 & 10 & {\bf 18} &    1 &       &         &         &            \\
 5             &   &   & 18 &   1 & {\bf 84} & {\bf 142} &       1 &         &            \\
 6             &   &   &  1 &  84 &  459 &     1 & {\bf 724} & {\bf 1266} &         1  \\
 7             &   &   &    & 142 &    1 &  5766 &   19057 &       1 &  {\bf 6516} \\
 8             &   &   &    &   1 &  724 & 19057 &       1 &  380597 &   1077681  \\
 9             &   &   &    &     & 1266 &     1 &  380597 & 3759277 &         1  \\
1      0       &   &   &    &     &    1 &  6516 & 1077681 &       1 & 185961668 
\end{tabular}
\end{center}
\end{table}

{\bf Problem 14:} Is it true that $T_{113} = 1$ when $a+b$ is 1 (mod 3)?
(Note that $a+b \equiv 1$ is the situation 
in which the Conway-Lagarias invariant
coincides with the area of the region being tiled,
so that all the tiles must be right-pointing stones.)

Comment: This problem is resolved in the affirmative 
by Theorem 1.1 of~\cite{DLPY}.


It is disappointing that the data for $a+b \not\equiv 1$ (mod 3)
do not suggest exact conjectures.
On the other hand, it is intriguing that congruence phenomena occur,
analogous to Cohn's 2-adic continuity theorem proved in~\cite{Cohn}
and conjectural 2-adic phenomena of a similar kind discussed in~\cite{Pro2}.

{\bf Problem 15:} Is $T_{113}(n,2n-4)$ 
2-adically continuous as a function of $n \geq 5$?

Comment: The sequence appears to be constant mod 2, constant mod 4,
2-periodic mod 8, and 8-periodic mod 16
(with repeating pattern 4, 4, 4, 4, 4, 12, 4, 12).
The case $n = 4$ breaks the pattern.
See \href{https://oeis.org/A364416}{OEIS entry A364416}.


{\bf Problem 16:} Is $T_{113}(n,2n-3)$ 
2-adically continuous as a function of $n \geq 4$?

Comment: The sequence appears to be constant mod 2, constant mod 4,
2-periodic mod 8, and 8-periodic mod 16
(with repeating pattern 2, 14, 2, 14, 10, 6, 10, 6).
The case $n = 3$ breaks the pattern.
See \href{https://oeis.org/A364417}{OEIS entry A364417}.


At the OPAC 2022 meeting, David Speyer suggested that 
one might use a different definition of a trimer,
namely, a path consisting of three vertices and two edges.
Thus, each stone would correspond to three different trimers
according to which 2 of the 3 possible edges one used.
Equivalently, one could use the original definition of a trimer,
with the proviso that each stone will count with weight 3.
The new numbers do not satisfy any nice patterns, with one exception:
in the case where all prototiles are allowed,
using stones of weight 3 seems to give rise to
the same 2-adic continuity phenomenon we saw in Problems 11 and 12.

Let $T_{ijk}(a,b;3)$ denote the weighted sum
of the $i,j,k$ tilings of the $(a,b)$-benzel
as defined near the end of section 1,
where a tiling with $m$ stones has weight $3^m$.


{\bf Problem 17:} Is $T_{113}(n,2n-4;3)$ 
2-adically continuous as a function of $n \geq 5$?

Comment: The sequence appears to be constant mod 2, constant mod 4,
constant mod 8, and 8-periodic mod 16
(with repeating pattern 4, 4, 12, 12, 4, 12, 12, 4).
The case $n = 4$ breaks the pattern.
See \href{https://oeis.org/A364418}{OEIS entry A364418}.


{\bf Problem 18:} Is $T_{113}(n,2n-3;3)$ 
2-adically continuous as a function of $n \geq 4$?

Comment: The sequence appears to be constant mod 2, constant mod 4,
2-periodic mod 8, and 8-periodic mod 16
(with repeating pattern 14, 6, 10, 2, 6, 14, 2, 10).
The case $n = 3$ breaks the pattern.
See 
\href{https://oeis.org/A364438}{OEIS entry A364438}.

Comment: Given that the number of right stones
minus the number of left stones in a tiling of a benzel
is the same for all tilings of that benzel, a different way
to describe the weighting assigned to tilings in Problems 17 and 18
(equivalent up to a power of 3)
is to give weight 1 to left stones and weight 9 to right stones.
During the final stages of preparing this article for publication,
I realized that one could more generally 
give weight 1 to left stones and weight $m$ to right stones.
The resulting enumerations appear to exhibit $p$-adic continuity
when $m+1$ is a power of the prime $p$.
Other forms of congruential consistency appear in the data as well;
for instance, the data for $m=9$ appear to satisfy 5-adic continuity.
I hope to address such patterns in future work.

\section{Miscellaneous} \label{sec-last}

\begin{figure}[h]
\begin{center}
\includegraphics[width=3.6in]{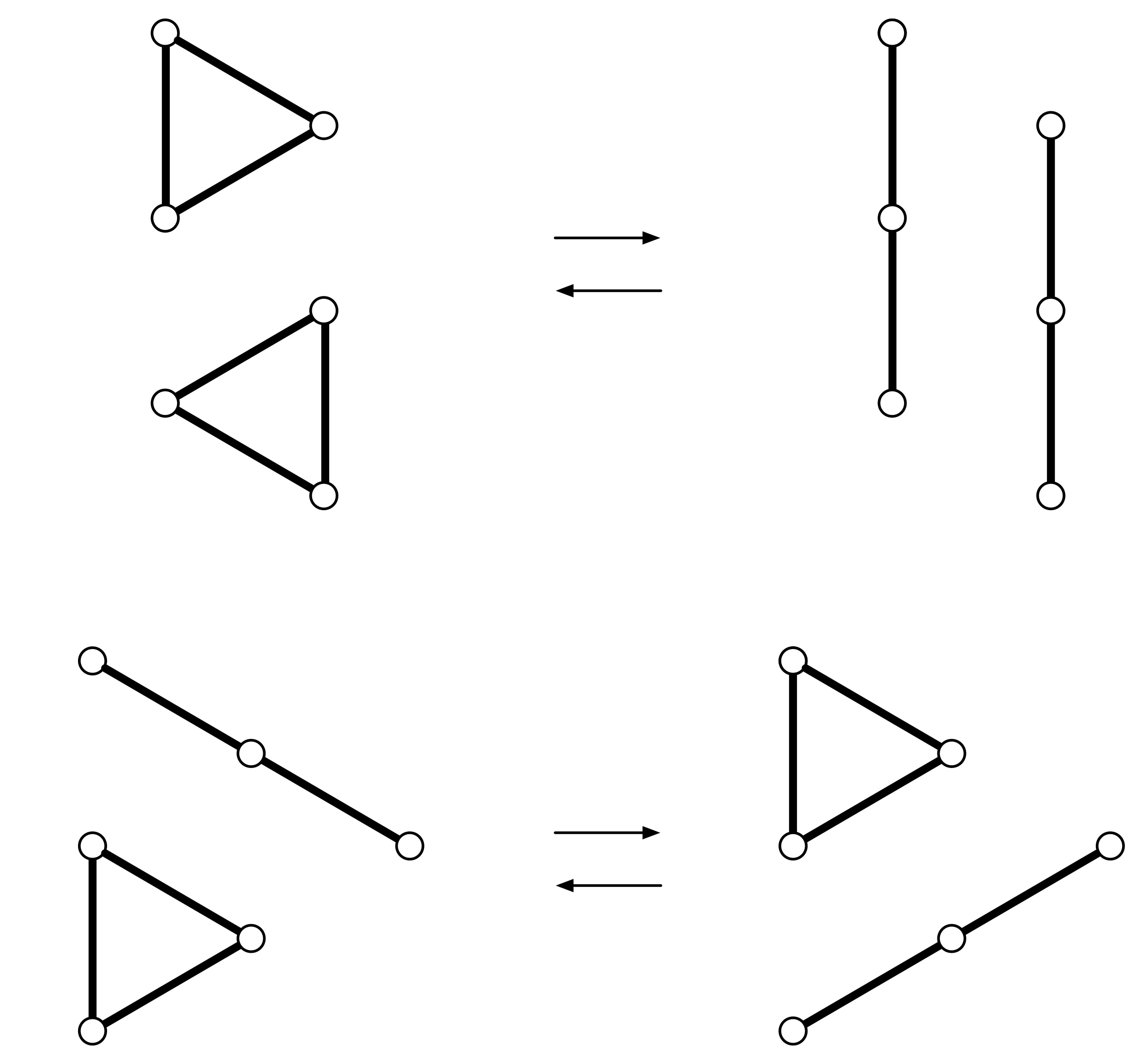}
\end{center}
\caption{2-flips for changing a trimer cover.}
\label{fig:mutate}
\end{figure}

The next question is not enumerative;
it is an old structural question that has gone unresolved for decades.
There are two natural kinds of ``2-flips''
that can turn one stones-and-bones tiling into another;
the first trades two stones of opposite orientation for two bones,
and the second trades a stone and a bone
for a stone of the same orientation
and a bone of a different orientation;
see Figure~\ref{fig:mutate}.

{\bf Problem 19:} Can every 
tiling of a finite simply-connected region
using stones and bones 
be mutated into every other such tiling
by means of a succession of 2-flips?

Comment: It is known that the hypothesis 
that the region be simply-connected cannot be dropped.
It is also known that if one restricts to tilings of type 112
(that is, if one prohibits one of the three orientations of bones),
then the claim is true; Sheffield proved an equivalent claim
in the context of ribbon tilings~\cite{Shef}.

\bigskip

In closing, we turn to the regions
that Conway, Lagarias, and Thurston originally studied
in papers~\cite{CoLa,Thur} that inspired much of my work on tilings:
triangles of hexagonal cells, with $n$ cells on each side (``$T_n$ regions'').
Those three authors showed that if one uses stones alone, 
$T_n$ can be tiled by $T_2$'s (that is, by stones)
precisely when $n$ is congruent to 0, 2, 9, or 11 (mod 12).
(In our notation, these are tilings of type 110;
such tilings were not discussed above since
for benzels they are not interesting from an enumerative perspective.)

The question we ask is, how many such tiling are there?
\href{https://oeis.org/A334875}{OEIS entry A334875} 
gives the answers for small values of $n$.
If we look at the prime factorizations of the answers,
we notice that the exponent of the prime 2 is creeping upward. 
Specifically, the multiplicity of the prime 2 
in the factorizations of the nonzero terms in this sequence 
goes 0, 0, 1, 3, 2, 3, 4, 3, 4, 3, 5, 8, 6, 8, \dots.
This is a priori surprising,
since the probability that a ``random'' positive integer
is divisible by $2^m$ decreases exponentially as $m$ increases.

{\bf Problem 20:} As $n$ goes to infinity
within the set of natural numbers 
congruent to 0, 2, 9, or 11 (mod 12),
does the number of tiling of $T_n$ by stones
converge 2-adically to 0?


\bigskip
\noindent
During the research described in this paper
I was supported by the Simons Foundation
Travel Support for Mathematicians program.
I am grateful to the organizers of OPAC 2022 
for inviting me to give this presentation.
I also thank Colin Defant, Rupert Li, and Ben Young
for useful conversations.

\bibliographystyle{amsplain}

\bibliographystyle{amsalpha}

\begin{thebibliography}{8}

\bibitem{Cohn}
Henry Cohn,
\href{https://www.combinatorics.org/ojs/index.php/eljc/article/view/v6i1r14}
{2-adic behavior of numbers of domino tilings},
Electron.\ J.\ Combin.\ {\bf 6} (1999), \#R14.

\bibitem{CoLa} 
John Conway and Jeffrey Lagarias, 
\href{http://dx.doi.org/10.1016/0097-3165(90)90057-4}
{Tiling with polyominoes and combinatorial group theory}, 
J.\ Combin.\ Theory (Ser.\ A) {\bf 53} (1990), 
no.\ 2, 183--208.

\bibitem{DLPY} 
Colin Defant, Rupert Li, James Propp, and Benjamin Young,
\href{https://escholarship.org/uc/item/5h16p4t1}
{Tilings of benzels via the abacus bijection},
Combinatorial Theory {\bf 3}(2) (2023).



\bibitem{KiPr} 
Jesse Kim and James Propp, 
\href{https://www.combinatorics.org/ojs/index.php/eljc/article/view/v30i3p26/pdf}
{A pentagonal number theorem for tribone tilings},
Electronic Journal of Combinatorics {\bf 30}(3) (2023), \#P3.26.

\bibitem{Pak} 
Igor Pak, 
\href{https://www.jstor.org/stable/221901}
{Ribbon tile invariants},
Trans.\ Amer.\ Math.\ Soc.\ {\bf 352} (2000), no.\ 12, 5525--5561.

\bibitem{Pro1}
James Propp,
\href{https://faculty.uml.edu/jpropp/eot.pdf}
{Enumeration of Tilings},
The Handbook of Enumerative Combinatorics (CRC Press, 2015), 541--588.

\bibitem{Pro2} 
James Propp, 
\href{http://math.colgate.edu/~integers/x30/x30.pdf}
{Some 2-adic conjectures concerning polyomino tilings of Aztec diamonds},
Integers {\bf 23} (1923), \#A30.

\bibitem{Shef} 
Scott Sheffield, 
\href{https://www.ams.org/journals/tran/2002-354-12/S0002-9947-02-02981-1/S0002-9947-02-02981-1.pdf}
{Ribbon tilings and multidimensional height functions},
Trans.\ Amer.\ Math.\ Soc.\ {\bf 354} (2002), no.\ 12, 4789--4813.


\bibitem{Thur} 
William Thurston, 
\href{https://doi.org/10.1080/00029890.1990.11995660}
{Conway’s tiling groups}, 
Amer.\ Math.\ Monthly {\bf 97} (1990), no.\ 8, 757--773. 

\bibitem{VeNi} 
Alain Verberkmoes and Bernard Nienhuis, 
\href{https://journals.aps.org/pre/abstract/10.1103/PhysRevE.63.066122}
{Bethe Ansatz solution of triangular trimers on the triangular lattice},
Phys.\ Rev.\ E {\bf 63} (2001), 066122.

\end{thebibliography}

\end{document}